\documentclass[11pt,doublespace]{article}
\begin{document}
 \oddsidemargin1.5cm




\newtheorem{thm}{Theorem}[section]
\newtheorem{lem}[thm]{Lemma}
\newtheorem{cor}[thm]{Corollary}
\newtheorem{ex}[thm]{Example}
\newtheorem{prop}[thm]{Proposition}
\newtheorem{remark}[thm]{Remark}
\newtheorem{coun}[thm]{Counterexample}
\newtheorem{defn}[thm]{Definition}
\newtheorem{conj}{Conjecture}
\newtheorem{problem}{Problem}
\newcommand\ack{\section*{Acknowledgement.}}

\newcommand{\etal}{{\it et al. }}

\newcommand{\bbP}{{\rm I\hspace{-0.8mm}P}}
\newcommand{\bbE}{{\rm I\hspace{-0.8mm}E}}
\newcommand{\bbF}{{\rm I\hspace{-0.8mm}F}}
\newcommand{\bbI}{{\rm I\hspace{-0.8mm}I}}
\newcommand{\bbR}{{\rm I\hspace{-0.8mm}R}}
\newcommand{\bbRp}{{\rm I\hspace{-0.8mm}R}_+}
\newcommand{\bbN}{{\rm I\hspace{-0.8mm}N}}
\newcommand{\bbC}{{\rm C\hspace{-2.2mm}|\hspace{1.2mm}}}
\newcommand{\bbD}{{\rm I\hspace{-0.8mm}D}}
\newcommand{\bbQ}{\bf Q}
\newcommand{\bbZ}{{\rm \rlap Z\kern 2.2pt Z}}
\newcommand{\bbK}{{\rm I\hspace{-0.8mm}K}}

\newcommand{\matP}{{\bbP}}
\newcommand{\mattildeP}{\tilde{\bbP}}
\newcommand{\matPN}[1]{{\bbP}_{#1}}
\newcommand{\matPP}[1]{{\bbP}_{#1}^0}
\newcommand{\matE}{{\bbE}}
\newcommand{\mattildeE}{\tilde{\bbE}}
\newcommand{\matEP}[1]{{\bbE}_{#1}^0}
\newcommand{\matF}{{\bbF}}
\newcommand{\matR}{{\bbR}}
\newcommand{\matRp}{{\bbRp}}
\newcommand{\matN}{{\bbN}}
\newcommand{\matZ}{{\bbZ}}
\newcommand{\matI}{{\bbI}}
\newcommand{\matK}{{\bbK}}
\newcommand{\matQ}{{\bbQ}}
\newcommand{\matC}{{\bbC}}
\newcommand{\matD}{{\bbD}}

\newcommand{\calL}{{\cal L}}
\newcommand{\calM}{{\cal M}}
\newcommand{\calN}{{\cal N}}
\newcommand{\calF}{{\cal F}}
\newcommand{\calG}{{\cal G}}
\newcommand{\calD}{{\cal D}}
\newcommand{\calB}{{\cal B}}
\newcommand{\calH}{{\cal H}}
\newcommand{\calI}{{\cal I}}
\newcommand{\calP}{{\cal P}}
\newcommand{\calQ}{{\cal Q}}
\newcommand{\calS}{{\cal S}}
\newcommand{\calT}{{\cal T}}
\newcommand{\calC}{{\cal C}}
\newcommand{\calK}{{\cal K}}
\newcommand{\calX}{{\cal X}}
\newcommand{\cals}{{\cal S}}
\newcommand{\calE}{{\cal E}}

\newcommand{\koniecmat}{\,}

\newcommand{\eqd}{\ =_{\rm d}\ }
\newcommand{\toto}{\leftrightarrow}
\newcommand{\eqdistr}{\stackrel{\rm d}{=}}
\newcommand{\as}{\stackrel{\rm a.s.}{=}}
\newcommand{\convdistr}{\stackrel{\rm d}{\rightarrow}}
\newcommand{\convweak}{{\Rightarrow}}
\newcommand{\convas}{\stackrel{\rm a.s.}{\rightarrow}}
\newcommand{\convfidi}{\stackrel{\rm fidi}{\rightarrow}}
\newcommand{\convprob}{\stackrel{p}{\rightarrow}}
\newcommand{\deff}{\stackrel{\rm def}{=}}
\newcommand{\bis}{{'}{'}}
\newcommand{\Cov}{{\rm Cov}}
\newcommand{\Var}{{\rm Var}}
\newcommand{\Exp}{{\rm E}}

\newcommand{\nd}{n^{\delta}}
\newcommand{\koniec}{\newline\vspace{3mm}\hfill $\odot$}


\title{Empirical process of long-range dependent sequences when parameters are estimated \protect}
\author{Rafa{\l} Kulik\thanks{School of Mathematics and Statistics,
University of Sydney, NSW 2006, Australia, email:
rkuli@maths.usyd.edu.au and Mathematical Institute, Wroc{\l}aw
University, Pl. Grunwaldzki 2/4, 50-384 Wroc{\l}aw, Poland}}

\maketitle

\begin{abstract}
In this paper we study the asymptotic behaviour of empirical
processes when parameters are estimated, assuming that the
underlying sequence of random variables is long-range dependent. We
show completely different phenomena compared to i.i.d. situation, as
well as compared to ordinary empirical processes of long range
dependent sequences. Applications include Kolmogorov-Smirnov and
Cramer-Smirnov-von Mises goodness-of-fit statistics.
\end{abstract}

\noindent {\bf Short title:} Estimated empirical processes and LRD

\section{Introduction and statement of results}
Let $\{\epsilon_i,i\ge 1\}$ be a centered sequence of i.i.d. random
variables. Consider the class of stationary linear processes
\begin{equation}\label{model}
X_i=\sum_{k=0}^{\infty}c_k\epsilon_{i-k},\ \ \ i\ge 1 .
\end{equation}
We assume that the sequence $c_k$, $k\ge 0$, is regularly varying
with index $-\beta$, $\beta\in (1/2,1)$ (written as $c_k\in
RV_{-\beta}$). This means that $c_k\sim k^{-\beta}L_0(k)$ as
$k\to\infty$, where $L_0$ is a slowly varying function at infinity.
We shall refer to all such models as long range dependent (LRD)
linear processes. In particular, if the variance exists, then the
covariances $\rho_k:=\Exp X_0X_k$ decay at the hyperbolic rate,
$\rho_k=L(k)k^{-(2\beta-1)}=:L(k)k^{-D}$, where
$\lim_{k\to\infty}L(k)/L_0^2(k)=B(2\beta-1,1-\beta)$ and
$B(\cdot,\cdot)$ is the beta-function. Consequently, the covariances
are not summable (cf. \cite{GiraitisSurgailis2002}).

Assume that $X_1$ has a continuous distribution function $F$. Given
$X_1,\ldots,X_n$, let $F_n(x)=n^{-1}\sum_{i=1}^n1_{\{X_i\le x\}}$ be
the empirical distribution function.

Assume that $E\epsilon_1^2<\infty$. Let $r$ be an integer and define
$$Y_{n,r}=\sum_{i=1}^n\sum_{1\le j_1<\cdots\le
j_r}\prod_{s=1}^rc_{j_s}\epsilon_{i-j_s},\qquad n\ge 1,$$ so that
$Y_{n,0}=n$, and $Y_{n,1}=\sum_{i=1}^nX_i$. If $p<(2\beta-1)^{-1}$,
then
\begin{equation}\label{eq-varance-behaviour}
\sigma_{n,p}^2:=\Var (Y_{n,p})\sim n^{2-p(2\beta-1)}L_0^{2p}(n).
\end{equation}
From \cite{HoHsing} we know that for $p<(2\beta-1)^{-1}$, as
$n\to\infty$,
\begin{equation}\label{eq-Ynr}
\sigma_{n,p}^{-1}Y_{n,p}\convdistr Z_p,
\end{equation}
where $Z_p$ is a random variable which can be represented by
appropriate multiple Wiener-It\^{o} integrals. In particular, $Z_1$
is standard normal.

In the present paper we study the asymptotic behaviour of empirical
processes when unknown parameters of the underlying distribution
function are estimated. The motivation to study such problems comes
from Kolmogorov-Smirnov type statistics. From \cite{HoHsing} we know
that, as $n\to\infty$,
\begin{equation}\label{eq-HoHsing}
\sigma_{n,1}^{-1}n\sup_{x\in {\matR}}|F_n(x)-F(x)|\convdistr
|Z_1|\sup_{x\in {\matR}}f(x),
\end{equation}
where $Z_1$ is a standard normal random variable and $f$ is the
density function of $F$. The above result can be used, in principle,
to test whether data $X_1,\ldots,X_n$ are consistent with a given
distribution $F$. If however $F$ belongs to a one-parameter family
$\{F(\cdot,{\theta}),\theta\in {\matR} \}$ say, then in order to use
(\ref{eq-HoHsing}) one needs to know the value of the parameter
$\theta$. A straightforward procedure would be to estimate it and
use the statistic
$$
\sigma_{n,1}^{-1}n\sup_{x\in {\matR}}|F_n(x)-F(x;\hat\theta_n)|,
$$
where $F(x;\hat\theta_n)$ is the distribution function
$F(x)=F(x;\theta)$ in which the parameter $\theta$ has been replaced
with its estimator $\hat\theta_n$. However, in the i.i.d. case, it
is known that such procedure changes a limiting process. To be more
specific, assume for a while that $X_1,\ldots,X_n$ are i.i.d. random
variables and consider
$$
\sqrt{n}\sup_{x\in {\matR}}|F_n(x)-F(x)|.
$$
As it is well-known, the above supremum converges in distribution to
the supremum of a Brownian bridge on $[0,1]$. On the other hand, for
a large class of estimators,
$$
\sqrt{n}|F_n(x)-F(x;\hat\theta_n)|,
$$
converges weakly to a Gaussian process, but no longer to a Brownian
bridge. The corresponding comments apply to the
Cram\'{e}r-Smirnov-von Mises statistic
$$
\sqrt{n}\int_{{\matR}}(F_n(x)-F(x))^2dF(x)
$$
and its 'estimated' version
$$
\sqrt{n}\int_{{\matR}}(F_n(x)-F(x;\hat\theta_n))^2dF(x;\hat\theta_n).
$$
We refer to \cite{darling}, \cite{durbin},
\cite{kac-kiefer-wolfowitz} and \cite{burke-csorgo-csorgo-revesz}
for more details.

Coming back to LRD sequences, we will focus on a location-scale
family of distributions. We shall assume that $Y_i=\sigma X_i+\mu$,
where $X_i$ is given by (\ref{model}) and $\sigma\not=0$. Clearly,
if $F$ is the distribution of $X_1$ and $H$ is the distribution of
$Y_1$, then $H(x)=F\left(\frac{x-\mu}{\sigma}\right)$. Moreover, the
empirical processes
$$
\beta_n(x)=\sigma_{n,1}^{-1}n(F_n(x)-F(x)),\qquad x\in {\matR}
$$
and
$$
\gamma_n(x)=\sigma_{n,1}^{-1}n(H_n(x)-H(x)),\qquad x\in {\matR}
$$
associated with $X_i$ and $Y_i$, respectively, are related by
\begin{equation}\label{eq-relation}
\gamma_n(x)=\beta_n\left(\frac{x-\mu}{\sigma}\right).
\end{equation}
From \cite{HoHsing}, $\beta_n(x)\convweak f(x)Z_1$, so that
$\gamma_n(x)\convweak f(\frac{x-\mu}{\sigma})Z_1$. Here and in the
sequel, $\convweak$ denotes weak convergence in
$D((-\infty,\infty))$. On the contrary, if $\hat\theta_n$ is an
appropriate sequence of estimators of the mean $\mu$, we will show
that, as $n\to\infty$,
$$\hat\gamma_n(x)=\sigma_{n,1}^{-1}n(H_n(x)-H(x;\hat\theta_n)),\qquad x\in {\matR}
$$
converges in probability to 0. Choosing a different scaling one can
obtain weak convergence, however the limiting process depends on the
choice of the estimator. In particular, using $\hat\theta_n=\bar
Y_n$ (the sample mean of $Y_1,\ldots,Y_n$) or $\hat\theta_n=M_n$
($M$-estimator), we can obtain different limits, depending on the
so-called {\it second-order M-rank} of the estimator $M_n$
introduced in \cite{KoulSurgailis}. Also, the scaling and the
limiting process depend on whether $\beta>3/4$ or $\beta<3/4$. In
particular, if $\beta>3/4$, then we obtain $\sqrt{n}$-consistency of
a modified Kolmogorov-Smirnov type statistics. The appropriate
results are stated in Theorems \ref{thm-beta-small} and
\ref{thm-beta-big}.

The proofs of our results will be based on a reduction principle for
long-range dependent empirical processes (see Theorem
\ref{thm-HoHsing} below), combined with approximation method as in
\cite{burke-csorgo-csorgo-revesz}. The fact, that we were able to
use the latter, Hungarian-like approach, shows its extreme power.
The Hungarian construction approach was for example employed to
obtain the Koml\'{o}s-Major-Tusn\'{a}dy (KMT) strong approximation
of empirical processes. Then, this approach was followed to
establish a number of optimal or almost optimal results for
functionals of empirical and quantile processes, including the one
in \cite{burke-csorgo-csorgo-revesz} for empirical processes with
parameters estimated (we refer to \cite{CsorgoHorvath}). The KMT
construction is tailored for the i.i.d. situation. However, a lot of
further developments based on this kind of approach, can be applied
to long-range dependent sequences. Very recent examples of such an
approach include \cite{CSW2006}, \cite{CsorgoKulik2006a},
\cite{kulik2006}.

The reduction principle was obtained first in
\cite{DehlingTaqqu1989} in case of subordinated Gaussian processes.
In more generality, it was obtained in the landmark paper
\cite{HoHsing}; see also \cite{KoulSurgailis2002} for related
studies. The best available result along these lines is due to Wu
\cite{Wu2003}. To state a particular version of his result, we shall
introduce the following assumptions, which will be valid throughout
the paper. Let $F_{\epsilon}$ be the distribution function of the
centered i.i.d. sequence $\{\epsilon_i,i\ge 1\}$. Assume that for a
given integer $p$, the derivatives $F^{(1)}_{\epsilon}, \ldots,
F^{(p+3)}_{\epsilon}$ of $F_{\epsilon}$ are bounded and integrable.
Note that these properties are inherited by the distribution $F$ as
well (cf. \cite{HoHsing} or \cite{Wu2003}).
\begin{thm}\label{thm-HoHsing}
Let $p$ be a positive integer. Then, as $n\to\infty$,
$$
\Exp \sup_{x\in {\matR}}\left|\sum_{i=1}^n(1_{\{X_i\le
x\}}-F(x))+\sum_{r=1}^p(-1)^{r-1}F^{(r)}(x)Y_{n,r}\right|^2=O(\Xi_n+n(\log
n)^2),
$$
where
$$
\Xi_n=\left\{\begin{array}{ll} O(n), & (p+1)(2\beta-1)>1\\
O(n^{2-(p+1)(2\beta-1)}L_0^{2(p+1)}(n)), & (p+1)(2\beta-1)<1
\end{array}\right. .
$$
\end{thm}
We will a require second-order expansion, thus in the above theorem,
$p=2$.

Let $\psi$ be a real-valued function of bounded variation such that
$\Exp{\psi(Y_1-\mu)}=0$. $M$-estimators are defined as
$$
M=M_n={\rm arg}\min \left\{\left|\sum_{j=1}^n\psi(Y_j-x)\right|,x\in
{\matR}\right\}.
$$
For $k=1,2$, let
$$
\lambda_k=\int_{\matR}\psi(y)f^{(k)}(y)dy .
$$
Let $k^*=k^*(\beta)=[1/(2\beta-1)]$, where $[\cdot]$ denotes the
integer part. The second-order rank $r_M(2)$ of the $M$-estimator
is: $r_M(2)=2$ if $k^*=1$ (so that $\beta>3/4$); $r_M(2)=2$ if
$k^*>1$ and $\lambda_2\not=0$; $r_M(2)>2$ if $k^*>1$ and
$\lambda_2=0$. We refer to \cite{KoulSurgailis} for more details. \\

Let
$$
a_n=\sigma_{n,2}\sigma_{n,1}^{-1}.
$$

Now, we are ready to state our results. We start with the case
$\beta<3/4$.
\begin{thm}\label{thm-beta-small}
Assume that $\theta_0=\mu$ and $\beta<3/4$. Then, under the
conditions of {\rm Theorem \ref{thm-HoHsing}}, as $n\to\infty$, we
have
\begin{itemize}
\item If $\hat\theta_n=\bar Y_n$ or $\hat\theta_n=M_n$, then
\begin{equation}\label{eq-negligible}
\sup_{x\in {\matR}}|\hat\gamma_n(x)|=o_P(1).
\end{equation}
\item If $\hat\theta_n=\bar Y_n$, then
\begin{equation}\label{eq-result-mean}
a_n^{-1}\hat\gamma_n(x)=\sigma_{n,2}^{-1}n(H_n(x)-H(x;\hat\theta_n))\convweak
f^{(1)}\left(\frac{x-\mu}{\sigma}\right)V,
\end{equation}
where $V$ is a linear combination of $Z_2$ and $\frac{1}{2}Z_1^2$.
\item If $\hat\theta_n=M_n$, $\Exp{\epsilon_1^{4\vee 2k^{*}(\theta)}}<\infty$ and $r_M(2)>2$, then
{\rm (\ref{eq-result-mean})} holds.
\item If $\hat\theta_n=M_n$, $\Exp{\epsilon_1^{4\vee 2k^{*}(\theta)}}<\infty$ and $r_M(2)=2$
\begin{equation}\label{eq-result-M}
a_n^{-1}\hat\gamma_n(x)=\sigma_{n,2}^{-1}n(H_n(x)-H(x;\hat\theta_n))\convweak
f^{(1)}\left(\frac{x-\mu}{\sigma}\right)V-\frac{\lambda_2}{2\lambda_1}\frac{1}{\sigma}f\left(\frac{x-\mu}{\sigma}\right)V_1,
\end{equation}
where $V$ is as in {\rm (\ref{eq-result-mean})} and $V_1$ is a
linear combination of $Z_1^2$ and $Z_2$.
\end{itemize}
\end{thm}
\begin{ex}{\rm
Assume that $\mu=0$, $f$ is symmetric and $\psi$ is skew-symmetric.
For $\beta<3/4$, $r_M(2)\ge 3$ (cf. \cite{KoulSurgailis}) and the
limiting behaviour is described by (\ref{eq-result-mean}). If,
however, $f$ is not symmetric, then $\lambda_2\not=0$ and
(\ref{eq-result-M}) holds. }
\end{ex}
As for the case $\beta>3/4$ we have the following theorem.
\begin{thm}\label{thm-beta-big}
Assume that $\theta_0=\mu$ and $\beta>3/4$. Then, under the
conditions of {\rm Theorem \ref{thm-HoHsing}}, as $n\to\infty$, we
have
\begin{itemize}
\item If $\hat\theta_n=\bar Y_n$ or $\hat\theta_n=M_n$, then
$$
\sup_{x\in {\matR}}|\hat\gamma_n(x)|=o_P(1).
$$
\item If $\hat\theta_n=\bar Y_n$, then
\begin{equation}\label{eq-result-mean-sqrtn}
\sqrt{n}\sigma_{n,1}n^{-1}\hat\gamma_n(x)=\sqrt{n}(H_n(x)-H(x;\hat\theta_n))\convweak
W\left(\frac{x-\mu}{\sigma}\right),
\end{equation}
where $W(\cdot)$ is a Gaussian process.
\item If $\hat\theta_n=M_n$, $\Exp{\epsilon_1^{4\vee
2k^{*}(\theta)}}<\infty$, then
\begin{equation}\label{eq-result-M-sqrtn}
\sqrt{n}\sigma_{n,1}n^{-1}\hat\gamma_n(x)=\sqrt{n}(H_n(x)-H(x;\hat\theta_n))\convweak
W\left(\frac{x-\mu}{\sigma}\right)+\frac{\sigma_{\psi}^2}{\sigma}f\left(\frac{x-\mu}{\sigma}\right)Z_1,
\end{equation}
$\sigma_{\psi}^2$ is given by the formula {\rm (1.18)} in {\rm
\cite{KoulSurgailis}}.
\end{itemize}
\end{thm}
An immediate corollary to Theorem \ref{thm-beta-small} is the
following Cram\'{e}r-Smirnov-von Mises test. An appropriate version
can also be stated in terms of Theorem \ref{thm-beta-big}.
\begin{cor}\label{cor-vonMises}
Let $\theta_0=\mu$ and $\hat\theta_n=\bar Y_n$. Under the conditions
of {\rm Theorem \ref{thm-beta-small}},
$$
\sigma_{n,2}^{-1}n\int_{\matR}(H_n(x)-H(x;\hat\theta_n))^2dH(x;\hat\theta_n)\convdistr
\frac{1}{\sigma}V^2\int_{{\matR}}\left(f^{(1)}\left(\frac{x-\mu}{\sigma}\right)\right)^2f\left(\frac{x-\mu}{\sigma}\right)dx.
$$
\end{cor}
The above result should be compared with a regular situation of non-estimated Cramer-Smirnov-von Mises statistics in
\cite{DehlingTaqqu1991}. The limiting distribution for the model
(\ref{model}) in case of Gaussian errors $\epsilon_i$, is a random
variable $Z_1^2$ multiplied by a deterministic function.\\

In what follows $C$ will denote a generic constant which may be
different at each of its appearance. Also, for any sequences $a_n$
and $b_n$, we write $a_n\sim b_n$ if $\lim_{n\to\infty}a_n/b_n=1$.
Moreover, $f^{(k)}$ denotes the $k$th order derivative of $f$.\\
\section{Proofs}
Let $p$ be a positive integer. Recall that
$$
a_n=\sigma_{n,2}\sigma_{n,1}^{-1}L_0(n),
$$
and let
$$
d_{n,p}=\left\{\begin{array}{ll} n^{-(1-\beta)}L_0^{-1}(n)(\log n)^{5/2}(\log\log n)^{3/4}, & (p+1)(2\beta-1)>1\\
n^{-p(\beta-\frac{1}{2})}L_0^{p}(n)(\log n)^{1/2}(\log\log n)^{3/4},
& (p+1)(2\beta-1)<1
\end{array}\right.
$$
Note that $d_{n,2}=o(a_n)$ provided $\beta<\frac{3}{4}$,

Put
\begin{eqnarray*}
S_{n,p}(x) &=&\sum_{i=1}^n(1_{\{X_i\le
x\}}-F(x))+\sum_{r=1}^p(-1)^{r-1}F^{(r)}(x)Y_{n,r}\\
&=:&\sum_{i=1}^n(1_{\{X_i\le x\}}-F(x))+V_{n,p}(x).
\end{eqnarray*}
Using Theorem \ref{thm-HoHsing} we obtain
\begin{eqnarray*}
\lefteqn{\sigma_{n,p}^{-1}\sup_{x\in
{\matR}}|S_{n,p}(x)|=}\\
&&\hspace*{-1cm}\left\{\begin{array}{ll} O_{a.s.}(n^{-(\frac{1}{2}-p(\beta-\frac{1}{2}))}L_0^{-p}(n)(\log n)^{5/2}(\log\log n)^{3/4}), & (p+1)(2\beta-1)>1\\
O_{a.s.}(n^{-(\beta-\frac{1}{2})}L_0(n)(\log n)^{1/2}(\log \log
n)^{3/4}), & (p+1)(2\beta-1)<1
\end{array}\right. .
\end{eqnarray*}
Since (see (\ref{eq-varance-behaviour}))
\begin{equation}\label{variances}
\frac{\sigma_{n,p}}{\sigma_{n,1}}\sim
n^{-(\beta-\frac{1}{2})(p-1)}L_0^{p-1}(n),
\end{equation}
we obtain
\begin{eqnarray}
\lefteqn{\sup_{x\in {\matR}}|\beta_n(x)+\sigma_{n,1}^{-1}V_{n,p}(x)|=}\label{eq-reduction-principle}\\
& = & \frac{\sigma_{n,p}}{\sigma_{n,1}}\sup_{x\in
{\matR}}\left|\sigma_{n,p}^{-1}\sum_{i=1}^n(1_{\{X_i\le
x\}}-F(x))+\sigma_{n,p}^{-1}V_{n,p}(x)\right|=o_{a.s.}(d_{n,p})\nonumber.
\end{eqnarray}

For a function $g(x;\theta)$ denote by
$\nabla_{\theta}^{r}g(x;\theta_0)$ its $r$th order derivative with
respect to $\theta$, evaluated at $\theta=\theta_0$. In particular,
$\nabla=\nabla^1$.

\subsection{Proof of Theorem \ref{thm-beta-small}}
Recall (\ref{eq-relation}). For an arbitrary unknown parameter
$\theta_0$ and its estimator $\hat\theta_n$ we have by
(\ref{eq-reduction-principle})
\begin{eqnarray}
\hat\gamma_n(x) & = &
\gamma_n(x)+\sigma_{n,1}^{-1}n(H(x;\theta_0)-H(x;\hat\theta_n))\nonumber\\
&= &
\beta_n\left(\frac{x-\mu}{\sigma}\right)+\sigma_{n,1}^{-1}n(H(x;\theta_0)-H(x;\hat\theta_n))\nonumber\\
&=
&o_p(d_{n,2})-\sigma_{n,1}^{-1}V_{n,2}\left(\frac{x-\mu}{\sigma}\right)+\sigma_{n,1}^{-1}n(\theta_0-\hat\theta_n)\nabla_{\theta}H(x;\theta_0)\nonumber\\
&&+\frac{1}{2}\sigma_{n,1}^{-1}n(\theta_0-\hat\theta_n)^2\nabla_{\theta}^2H(x;\theta_0)+
\frac{1}{6}\sigma_{n,1}^{-1}n(\theta_0-\hat\theta_n)^3\nabla_{\theta}^3H(x;\hat\theta_n^*)\nonumber\\
&= &
o_p(d_{n,2})-\sigma_{n,1}^{-1}f\left(\frac{x-\mu}{\sigma}\right)\sum_{i=1}^nX_i+\sigma_{n,1}^{-1}f^{(1)}\left(\frac{x-\mu}{\sigma}\right)Y_{n,2}\nonumber\\
&&+\sigma_{n,1}^{-1}n(\theta_0-\hat\theta_n)\nabla_{\theta}H(x;\theta_0)+\frac{1}{2}\sigma_{n,1}^{-1}n(\theta_0-\hat\theta_n)^2\nabla_{\theta}^2H(x;\theta_0)\nonumber\\
&&+
\frac{1}{6}\sigma_{n,1}^{-1}n(\theta_0-\hat\theta_n)^3\nabla_{\theta}^3H(x;\hat\theta_n^*)\label{eq-technical-1},
\end{eqnarray}
with some $\hat\theta_n^*$ such that
$|\hat\theta_n^*-\hat\theta_n|\le |\theta_0-\hat\theta_n^*|$. 

If $\theta_0=\mu$, then
\begin{equation}\label{eq-derivative-mu}
\nabla_{\theta}^rH(x)=\nabla_{\mu}^rF\left(\frac{x-\mu}{\sigma}\right)=(-1)^r\frac{1}{\sigma^r}f^{(r-1)}\left(\frac{x-\mu}{\sigma}\right).
\end{equation}
Also, if $\hat\theta_n=\bar Y_n$, then
\begin{equation}\label{eq-estim-mean}
\hat\theta_n-\theta_0=\sigma\bar X_n
\end{equation}
Hence, using uniform boundness of $f^{(2)}$,
\begin{eqnarray*}
\hat\gamma_n(x) & = &
o_p(d_{n,2})-\sigma_{n,1}^{-1}f\left(\frac{x-\mu}{\sigma}\right)\sum_{i=1}^nX_i+\sigma_{n,1}^{-1}f^{(1)}\left(\frac{x-\mu}{\sigma}\right)Y_{n,2}+\\
&&\sigma_{n,1}^{-1}f\left(\frac{x-\mu}{\sigma}\right)\sum_{i=1}^nX_i+\frac{1}{2}\sigma_{n,1}^{-1}nf^{(1)}\left(\frac{x-\mu}{\sigma}\right)\bar
X_n^2+ O_P\left(\sigma_{n,1}^{-1}n\bar X_n^3 \right).
\end{eqnarray*}
Since $\beta<3/4$, note that $\sigma_{n,1}Y_{n,2}=o_p(1)$ (cf. (\ref{eq-Ynr})), $\sigma_{n,1}^{-1}n\bar X_n^2=o_P(1)$ and $\sigma_{n,1}^{-1}n\bar X_n^3=o_P(1)$.
Thus, we conclude that $\sup_{x}|\hat\gamma_n(x)|\convprob 0$ for $\hat\theta_n=\bar Y_n$. \\

Further,
\begin{eqnarray*}
\lefteqn{a_n^{-1}\sup_x\left|\hat\gamma_n(x)-f^{(1)}\left(\frac{x-\mu}{\sigma}\right)\left[\sigma_{n,1}^{-1}Y_{n,2}+\frac{1}{2}\sigma_{n,1}^{-1}n\bar
X_n^2\right]\right|}\\
&=&o_p(d_{n,2}a_n^{-1})+O_P(a_n^{-1}\sigma_{n,1}^{-1}n\bar X_n^3) =
o_p(1)+O_P(a_n^{-1}\sigma_{n,1}^{-1}nn^{-3}\sigma_{n,1}^3)\\
&=&o_P(1).
\end{eqnarray*}
Thus, (\ref{eq-result-mean}) follows.\\

If $\hat\theta_n=M_n$ then, as in (\ref{eq-technical-1}) and
(\ref{eq-derivative-mu}),
\begin{eqnarray*}
\hat\gamma_n(x) & = & o_p(d_{n,2})-\sigma_{n,1}^{-1}f\left(\frac{x-\mu}{\sigma}\right)\sum_{i=1}^nX_i+\sigma_{n,1}^{-1}f^{(1)}\left(\frac{x-\mu}{\sigma}\right)Y_{n,2}+\\
&& -\frac{1}{\sigma}\sigma_{n,1}^{-1}n(\mu-\bar
Y_n)f\left(\frac{x-\mu}{\sigma}\right)-\frac{1}{\sigma}\sigma_{n,1}^{-1}n(\bar
Y_n-M_n)f\left(\frac{x-\mu}{\sigma}\right)+\\
&&\frac{1}{2\sigma^2}\sigma_{n,1}^{-1}nf^{(1)}\left(\frac{x-\mu}{\sigma}\right)(\mu-M_n)^2+O_P(\sigma_{n,1}^{-1}n(\mu-M_n)^3)\\
&= &
o_p(d_{n,2})+\sigma_{n,1}^{-1}f^{(1)}\left(\frac{x-\mu}{\sigma}\right)Y_{n,2}-\frac{1}{\sigma}\sigma_{n,1}^{-1}n(\bar
Y_n-M_n)f\left(\frac{x-\mu}{\sigma}\right)\\
&&+\frac{1}{2\sigma^2}\sigma_{n,1}^{-1}nf^{(1)}\left(\frac{x-\mu}{\sigma}\right)(\mu-M_n)^2+O_P(\sigma_{n,1}^{-1}n(\mu-M_n)^3).
\end{eqnarray*}
From \cite{KoulSurgailis},
\begin{equation}\label{eq-conv-M}
\sigma_{n,1}^{-1}n(M_n-\mu)=\sigma_{n,1}^{-1}n(\bar
Y_n-\mu)+o_P(1)\convdistr \sigma^2Z_1
\end{equation}
and $\sigma_{n,1}^{-1}n(\bar Y_n-M_n)=o_P(1)$. Thus,
$\sup_{x}|\hat\gamma_n(x)|\convprob 0$ for $\hat\theta_n=M_n$.\\

If $r_M(2)>2$, then from \cite[Theorem 1.1]{KoulSurgailis},
$$
a_n^{-1}\sigma_{n,1}^{-1}n(\bar Y_n-M_n)=o_P(1),
$$
thus in this case
\begin{eqnarray*}
\lefteqn{a_n^{-1}\sup_x\left|\hat\gamma_n(x)-f^{(1)}\left(\frac{x-\mu}{\sigma}\right)\left[\sigma_{n,1}^{-1}Y_{n,2}+\frac{1}{2\sigma^2}\sigma_{n,1}^{-1}n(\mu-M_n)^2\right]\right|}\\
&=&o_p(d_{n,2}a_n^{-1})+o_P(1)+O_P(a_n^{-1}\sigma_{n,1}^{-1}n(\mu-M_n)^3)
=o_P(1).
\end{eqnarray*}
Therefore, in view of (\ref{eq-conv-M}), (\ref{eq-result-mean}) follows.

If $r_M(2)=2$, then $a_n^{-1}\sigma_{n,1}^{-1}n$ is the proper
scaling for $(\bar Y_n-M_n)$ and thus
\begin{eqnarray*}
\lefteqn{a_n^{-1}\sup_x\left|\hat\gamma_n(x)-f^{(1)}\left(\frac{x-\mu}{\sigma}\right)\left[\sigma_{n,1}^{-1}Y_{n,2}+\frac{n(\mu-M_n)^2}{2\sigma^2\sigma_{n,1}}\right]\right.}\\
&&\qquad\qquad
+\left.\frac{n}{\sigma\sigma_{n,1}}f\left(\frac{x-\mu}{\sigma}\right)(\bar
Y_n-M_n)\right|\\
&=&o_p(d_{n,2}a_n^{-1})+O_P(a_n^{-1}\sigma_{n,1}^{-1}n(\mu-M_n)^3)
=o_P(1),\qquad\qquad\qquad\qquad \qquad\qquad
\end{eqnarray*}
and hence (\ref{eq-result-M}) follows using (\ref{eq-conv-M}) and Corollary 1.1 in \cite{KoulSurgailis}.\koniec
\subsection{Proof of Corollary \ref{cor-vonMises}}
Write
\begin{eqnarray*}
\lefteqn{\int\hat\gamma_n(x)^2dH(x;\hat\theta_n)=\int\hat\gamma_n(x)^2h(x;\theta_0)dx}\\
&&+\int\hat\gamma_n(x)^2(h(x;\hat\theta_n-h(x;\theta_0))dx .
\end{eqnarray*}
As for the second term, we have
$$
\int\hat\gamma_n(x)^2\nabla_{\theta}h(x;\theta_0)
(\hat\theta_n)-\theta_0)dx +R_n,
$$
where $R_n=O_P((\hat\theta_n-\theta_0)^2)=o_P(\hat\theta_n-\theta_0)$. Thus, the second term is of
a smaller rate than the first one and the limiting behaviour of
$a_n^{-1}\int\hat\gamma_n(x)^2dH(x;\hat\theta_n)$ is the same as
that of $\int\hat\gamma_n(x)^2h(x;\theta_0)dx$. Thus, Corollary
\ref{cor-vonMises} follows from Theorem \ref{thm-beta-small}.\koniec

\subsection{Proof of Theorem \ref{thm-beta-big}}
Recall that $\beta>3/4$. Then
\begin{eqnarray*}
\lefteqn{\sqrt{n}\sigma_{n,1}n^{-1}\hat\gamma_n(x)  =
\sqrt{n}\sigma_{n,1}n^{-1}\beta_n\left(\frac{x-\mu}{\sigma}\right)+\sqrt{n}\left(F\left(\frac{x-\mu}{\sigma}\right)-F\left(\frac{x-\mu}{\sigma},\hat\theta_n\right)\right)}\\
&= &
\sqrt{n}\left(F_n\left(\frac{x-\mu}{\sigma}\right)-F\left(\frac{x-\mu}{\sigma}\right)+f\left(\frac{x-\mu}{\sigma}\right)\sum_{i=1}^nX_i/n\right)\\
&&-f\left(\frac{x-\mu}{\sigma}\right)\frac{\sum_{i=1}^nX_i}{\sqrt{n}}-\frac{1}{\sigma}\sqrt{n}(\theta_0-\hat\theta_n)f\left(\frac{x-\mu}{\sigma}\right)+O(\sqrt{n}(\theta_0-\hat\theta_n)^2)\\
&:= & W_n\left(\frac{x-\mu}{\sigma}\right)-f\left(\frac{x-\mu}{\sigma}\right)\frac{\sum_{i=1}^nX_i}{\sqrt{n}}-\frac{1}{\sigma}\sqrt{n}(\theta_0-\hat\theta_n)f\left(\frac{x-\mu}{\sigma}\right)\\
&&+O(\sqrt{n}(\theta_0-\hat\theta_n)^2).
\end{eqnarray*}
If $\theta_0=\mu$ and $\hat\theta_n=\bar Y_n$, then via
(\ref{eq-estim-mean}),
$$
\sup_{x\in{\matR}}\left|\sqrt{n}\sigma_{n,1}n^{-1}\hat\gamma_n(x)-W_n\left(\frac{x-\mu}{\sigma}\right)\right|=O_P(\sqrt{n}(\mu-\hat\theta_n)^2)=o_P(1).
$$
Thus, using \cite[Theorem 3]{Wu2003}, we obtain
(\ref{eq-result-mean-sqrtn}).

If $\theta_0=\mu$ and $\hat\theta_n=M_n$, then
$$
\sup_{x\in{\matR}}\left|\sqrt{n}\sigma_{n,1}n^{-1}\hat\gamma_n(x)-W_n(x)+\frac{1}{\sigma}f\left(\frac{x-\mu}{\sigma}\right)\sqrt{n}(M_n-\bar
Y_n)\right| =o_P(1).
$$
If $\beta>3/4$, then from \cite[Theorem 1.1]{KoulSurgailis},
$\sqrt{n}(M_n-\bar Y_n)\convdistr N(0,\sigma_{\phi}^2)$. Thus,
(\ref{eq-result-M-sqrtn}) follows.\koniec

\ack This work was initiated during my stay at Carleton University.
I am thankful to Professors Barbara Szyszkowicz and Mikl\'{o}s
Cs\"{o}rg\H{o} for their support and helpful remarks.


\begin{thebibliography}{99}
\bibitem{burke-csorgo-csorgo-revesz}
Burke, M. D., Cs\"{o}rg\H{o}, M., Cs\"{o}rg\H{o}, S.,
R\'{e}v\'{e}sz, P. (1979). Approximations of the empirical process
when parameters are estimated.  {\it Ann. Probab.} {\bf 7},
790--810.
\bibitem{CsorgoHorvath}
Cs\"{o}rg\H{o}, M., Horv\'{a}th, L. (1993). {\it Weighted
approximations in probability and statistics}. Wiley Series in
Probability and Mathematical Statistics: Probability and
Mathematical Statistics. John Wiley \& Sons, Ltd., Chichester, 1993.
\bibitem{CSW2006}
Cs\"{o}rg\H{o}, M., Szyszkowicz, B. and Wang, L. (2006). Strong
Invariance Principles for Sequential Bahadur-Kiefer and Vervaat
Error Processes of Long-Range Dependence Sequences. {\it Ann.
Statist.} {\bf 34}, 1013--1044.
\bibitem{CsorgoKulik2006a}
Cs\"{o}rg\H{o}, M. and Kulik, R. (2006). Reduction principles for
quantile and Bahadur-Kiefer processes of long-range dependent linear
sequences. {\it Preprint}.
\bibitem{darling}
Darling, D. A. (1955). The Cram\'{e}r-Smirnov test in the parametric
case. {\it Ann. Math. Statist.} {\bf 26}, 1--20.
\bibitem{DehlingTaqqu1989}
Dehling, H. and Taqqu, M. (1989). The Empirical Process of some
Long-Range Dependent Sequences with an Applications to
$U$-Statisitcs. {\it Ann. Statist.} {\bf 17}, 1767--1783.
\bibitem{DehlingTaqqu1991}
Dehling, H. and Taqqu, M. (1991). Bivariate symmetric statistics of
long-range dependent observations. {\it J. Statist. Pl. Inf.} {\bf
28}, 153--165.
\bibitem{durbin}
Durbin, J. (1973). Weak convergence of the sample distribution
function when parameters are estimated.  {\it Ann. Statist.} {\bf
1}, 279--290.
\bibitem{GiraitisSurgailis2002}
Giraitis, L. and Surgailis, D. (2002). The reduction principle for
the empirical process of a long memory linear process. {\it
Empirical process techniques for dependent data}, 241--255,
Birkh\"{a}user Boston, Boston, MA.
\bibitem{HoHsing}
Ho, H.-C. and Hsing, T. (1996). On the asymptotic expansion of the
empirical process of long-memory moving averages. {\it Ann.
Statist.} {\bf 24}, 992--1024.
\bibitem{kac-kiefer-wolfowitz}
Kac, M., Kiefer, J., Wolfowitz, J. (1955). On tests of normality and
other tests of goodness of fit based on distance methods. {\it Ann.
Math. Statist.} {\bf 26}, 189--211.
\bibitem{KoulSurgailis}
Koul, H.L. and Surgailis, D. (1996). Asymptotic expansion of
$M$-estimators with long memory moving errors. {\it Ann. Statist.}
{\bf 25}, 818--850.
\bibitem{KoulSurgailis2002}
Koul, H.L. and Surgailis, D. (2002). Asymptotic expansion of the
empirical process of long memory moving averages. {\it Empirical
process techniques for dependent data}, 213--239, Birkh\"{a}user
Boston, Boston, MA.
\bibitem{kulik2006}
Kulik, R. (2006). Sums of extreme values of subordinated long-range
dependent sequences: moving averages with finite variance. {\it
Submitted.}
\bibitem{Wu2003}
Wu, W.B. (2003). Empirical processes of long-memory sequences. {\it
Bernoulli} {\bf 9}, 809--831.
\end{thebibliography}
\end{document}